\documentclass[11pt,twoside]{article}
\usepackage{amsfonts,amsmath,amssymb,amsthm}
\usepackage{epsfig}
\everymath{\displaystyle}
\def\l{\left}
\def\r{\right}

\newcommand{\bb}[1]{\begin{equation}\label{#1}}
\newcommand{\ee}{\end{equation}}
\newcommand{\bbb}{\begin{eqnarray}}
\newcommand{\eee}{\end{eqnarray}}
\newcommand{\bbbb}{\begin{eqnarray*}}
\newcommand{\eeee}{\end{eqnarray*}}
\newcommand{\nnn}{\nonumber}
\newcommand{\no}{\noindent}

\def\CC{\mathbb{C}}
\def\ZZ{\mathbb{Z}}
\def\R#1{$(\ref{#1})$}

\newtheorem{Theorem}{Theorem}[section]

\newtheorem{Corollary}[Theorem]{Corollary}

\newtheorem{Definition}[Theorem]{Definition}

\begin{document}
\baselineskip=4.8mm

\begin{center}

{\Large\bf An Effective $z$-Stretching Method for Paraxial Light
Beam Propagation Simulations\footnote{The third and fourth authors
are supported in part by a research grant (No. AFGD-035-75CS) from
the Air Force Research Laboratory and General Dynamics. The fourth
author is also supported in part by ASEE-SFFP Awards from the U.S.
Air Force.}}

\vspace{8mm}

$\mbox{Leonel Gonzalez}^{\dag},~\mbox{Shekhar Guha}^{\dag},~
\mbox{James W. Rogers}\footnote{Principal and corresponding
author. Email address:
James\underline{~}W\underline{~}Rogers@baylor.edu}^{\natural}~\mbox{and
Qin Sheng}^{\flat,\natural}$

\vspace{3mm}

$\mbox{~}^{\dag}$Materials \& Manufacturing Directorate, Air Force
Research Laboratory\\2977 Hobson Way, Wright-Patterson AFB, OH
45433-7702, USA

\vspace{2mm}

$\mbox{~}^{\flat}$Center for Astrophysics, Space Physics \& Engineering Research\\
$\mbox{~}^{\natural}$Department of Mathematics\\
Baylor University, Waco, TX 76798-7328, USA

\vspace{6mm}

\parbox[t]{13cm}{\small \noindent{\small{\bf Abstract.}
A $z$-stretching finite difference method is established for
simulating the paraxial light beam propagation through a lens in a
cylindrically symmetric domain. By introducing proper domain
transformations, we solve corresponding difference approximations on
a uniform grid in the computational space for great efficiency. A
specialized matrix analysis method is constructed to study the
numerical stability. Interesting computational results are
presented. }

\vspace{4mm}

\parbox[t]{13cm}{\small{\bf Keywords.} Light beam propagation, interface or discontinuous
surface, domain transformation, consistency, stability, uniform and
nonuniform grids, approximations}

\vspace{2mm}

\parbox[t]{13cm}{\small{\bf AMS (MOS) Subject Classification:} 65M06, 65M50, 65Z05, 78A15,
78M20}}
\end{center}

\vspace{8mm}

\noindent{\large\bf 1.~~Introduction}
\setcounter{section}{1}\setcounter{equation}{0}

\vspace{1mm}

In order to reduce the computational complexity of light beam
propagation simulation, a number of approximations are typically
employed.  We take advantage of these approximation techniques for
the paraxial case to derive an efficient and robust method that
allows the application of conventional finite difference schemes on
a uniform grid in the computational space, despite the difficulty of
an interface present in the domain.

From Maxwell's field equations describing the behavior of monochromatic light,
we obtain the time-dependent Helmholtz equation,
\bb{a1} \nabla^{2}E - \frac{1}{c^{2}}\frac{\partial^{2}E}{\partial
t^{2}}=0, \ee where $E=E(x,y,z,t)$ is the electric field
intensity, $\nabla^2$ is the Laplacian operator, and $c$ is the
phase velocity, or speed of light in a particular medium.

Let $E$ be the field intensity of a monochromatic plane wave of
the form \bb{ComplexWavefunction}E(x,y,z,t) = U(x,y,z)e^{i2\pi\nu
t},\ee where $\nu$ is the frequency of the light. Then from \R{a1}
we acquire the time-independent Helmholtz equation
\bb{Helmholtz}\left(\nabla^2 + \kappa^{2}\right)U(x,y,z) = 0,\ee
where $\kappa = {2\pi\nu}/{c} = {2\pi}/{\lambda}$ is referred as
the {\em wave number,} and $\lambda = {c}/{\nu}$ is referred as
the {\em wavelength.} Functions $|U(x,y,z)|$ and $\arg(U(x,y,z))$
are the {\em amplitude\/} and {\em phase of the wave,\/}
respectively. Assume that the $z$ is the direction of the beam
propagation. We may further consider the wave function with a
complex amplitude, that is, \bb{ModulatedPlaneWave}U(x,y,z) =
u(x,y,z)e^{-i\kappa z}\ee where $u$ is called a {\em complex
envelope.\/} A paraxial wave becomes realistic if the variation of
$u$ is slow in the $z$-direction.

Substitute \R{ModulatedPlaneWave} into \R{Helmholtz} to yield
\bb{AlmostParaxialHelmholtz}
\nabla^{2}_{T}u(x,y,z) - 2ik\frac{\partial u(x,y,z)}{\partial z} +
\frac{\partial^{2} u(x,y,z)}{\partial z^{2}} = 0\ee
where
$$\nabla^{2}_{T} = \frac{\partial^{2}}{\partial x^{2}} +
\frac{\partial^{2}}{\partial y^{2}}$$ is the transverse Laplacian
operator. In a paraxial case, we may assume that within a
wavelength of the propagation distance, the change in $u$ is
sufficiently small compared to $|u|$ \cite{Poon}. Thus,
$$\left|\frac{\partial^{2} u}{\partial z^{2}}\right|\ll
|\kappa^{2}u|$$ which indicates that
$$\frac{\partial^{2} u}{\partial z^{2}}\approx 0.$$
Therefore we arrive at an approximation of
\R{AlmostParaxialHelmholtz},
\bb{ParaxialHelmholtz}\nabla^{2}_{T}u(x,y,z) - 2i\kappa\frac{\partial
u(x,y,z)}{\partial z} = 0,\ee which is called the {\em slowly
varying envelope approximation of the Helmholtz equation\/} \cite{Band,Guha1,Saleh1}.

Under the transformation $r=\sqrt{x^{2} + y^{2}}$ and $\phi =
\arctan({y}/{x}),$ \R{ParaxialHelmholtz} can be reformulated to
\bb{PolarParaxialHelmholtz}
\left(\frac{1}{r}\frac{\partial}{\partial r} +
\frac{\partial^{2}}{\partial r^{2}} +
\frac{1}{r^{2}}\frac{\partial^{2}}{\partial \phi^{2}} -
2i\kappa\frac{\partial}{\partial z}\right)u(r,z) = 0. \ee The
equation (\ref{PolarParaxialHelmholtz}) has been utilized
frequently in laser beam propagation simulations in the past
decades \cite{Guha1,Guha2,Shib1}.

\begin{center}
\fbox{\epsfig{file=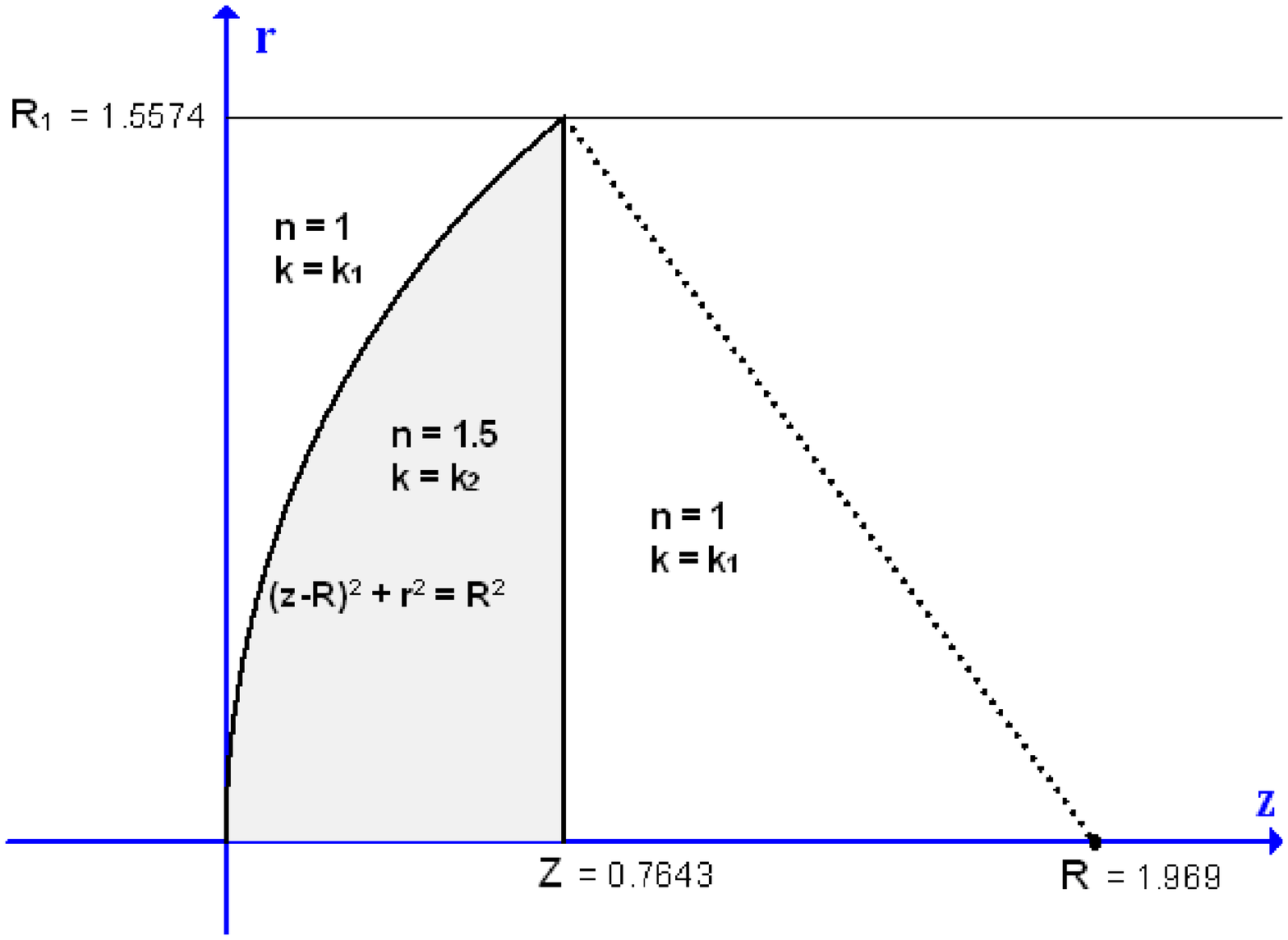,width=4.2in,height=2.5in}}
%\fbox{\epsfig{file=f060806_c01_lens2b.eps,width=4.2in,height=1.8in}}

\vspace{2mm}

\parbox[t]{13cm}{\small{\sc Figure 1.1}. ~An illustration of the
leans area \cite{Sheng2}. }
\end{center}

In this paper, we shall consider a cylindrically symmetric domain
for spherical lens environments. For the situation we may assume
\cite{Guha1} that
$$\frac{\partial^{2}u}{\partial \phi^{2}} \equiv 0$$
and \R{PolarParaxialHelmholtz} can be simplified to yield \bb{1f}
2i\kappa\frac{\partial u}{\partial z}(r,z) = \frac{\partial^2
u}{\partial r^2}(r,z) + \frac{1}{r}\frac{\partial u}{\partial
r}(r,z),~~~0 \leq r \leq r_{0} \ll \infty. \ee We assume that the
wavelength of light, $\lambda,$ is 9.449 $\mu m$ and that light is
incident from air (refractive index $n_1=1$) into glen (refractive
index $n_2=1.5$) so that we may adopt the following wave numbers
\cite{Guha1} \bb{1g} \kappa(r,z) = \left\{
\begin{array}{ll} {\kappa}_0 =\dfrac{2\pi n_{1}}{\lambda} \approx
\dfrac{2}{3}\times 9.97543 \times 10^{3} \mbox{\em cm}^{-1},& \mbox{ in medium one;}\\
{\kappa}_1 = \dfrac{2\pi n_{2}}{\lambda}\approx 9.97543 \times
10^{3}\mbox{\em cm}^{-1}, & \mbox{ in medium two.}
\end{array}\right.
\ee The above implies that a coefficient in
\R{PolarParaxialHelmholtz} is discontinuous at the lens interface.
Relation \R{1g} represents a single surface situation. Multiple
surface scenarios can also be discussed with additional $\kappa$
values. Needless to say, the discontinuity adds considerable
difficulties to the computation of the numerical solution of the
differential equation \cite{Cheng,Goodman,Larsson,Sheng2}.

%% CHANGE-2A replaced paragraph

%As for the initial condition for \R{PolarParaxialHelmholtz} and
%\R{1f}, based on the assumption of a complex electric field function
%of the Gaussian beam, we may assume that \bb{initial}u(r,0) =
%e^{-(r/\omega_{0})^{2}},\ee where $\omega_{0}$ is the Gaussian beam
%width. Further, based on the symmetry of the solution at $r=0$ and
%the fact that the optical waves do not exist at the top of the
%domain, we may introduce the following Neumann boundary conditions
%\bb{bound_cond}\frac{\partial u}{\partial r}(0,z)= \frac{\partial
%u}{\partial r}(r_{0},z)=0,~~~z>0. \ee

We employ Neumann boundary conditions
\begin{equation}\label{bound_cond}u_{r}(z,0)=u_{r}(z,R_{1})=0, ~~ z>0,
\end{equation}
\noindent at the bottom, $r=0$ and top, $r
= R_{1}$ of the rectangular domain.

For the initial solution of boundary value problem \R{PolarParaxialHelmholtz}, \R{1f}, we use the following approximation of a Gaussian beam with point source \cite{Eom}
\begin{equation}\label{init_val}u(z,r) = \frac{A}{1 + i\vartheta}\exp\left(ikz - \frac{r^{2}}{\beta^{2}(1+i\vartheta)}\right),
\end{equation}
\noindent where $\beta_{0}$ is the Gaussian beam
width, while $\vartheta$, $\beta$ and $A$ are parameters such that
$$\vartheta = \frac{2z}{\beta^{2}k}, ~~~\frac{1}{\beta^{2}} = \frac{1}{\beta_{0}^{2}} + \frac{ik}{2z_{0}}, ~~~ A = e^{ikz_{0}}.$$

\vspace{8mm}

\noindent{\large\bf 2.~~Base Difference Scheme and Stability}
\setcounter{section}{2}\setcounter{equation}{0}

\vspace{1mm}

%% CHANGE-2A replacing all r_0 with R_{1}, bound_cond with init_val

\noindent Let $0\leq z\leq Z$ for \R{1f}-\R{init_val}. Further,
let $h=R_{1}/M,~\tau=Z/N,$ where $M,~N \in \ZZ^{+}$ are sufficiently
large. We may introduce the uniform grid region,
$$\Omega_{h,\tau} = \left\{(mh,n\tau)~|~0 \leq m \leq M, 0 \leq n \leq N \right\}$$
over the rectangular domain $\Omega$ used. For the sake of
simplicity, we denote $r_{m} = mh$ and $z_{n} = n\tau.$ In
addition, we will use $z_{n-\alpha},~0<\alpha<1,$ for specifying a
non-grid point between $z_{n-1}$ and $z_{n}$ whenever needed.

%% CHANGE-2A bound_cond, init_val

Let us start with a linear second-order partial differential
equation of the form \bb{general_PDE} c_{5}\frac{\partial^2
u}{\partial z \partial r} + c_{4}\frac{\partial^2
u}{\partial r^2} + c_{3}\frac{\partial u}{\partial r} +
c_{2}\frac{\partial u}{\partial z} + c_{1}u + c_{0} =
0,~~~(r,z)\in\Omega, \ee together with \R{bound_cond},
\R{init_val}.  The coefficients $c_{i}$ of \R{general_PDE} are functions of $z$ and $r$ and may
be discontinuous due to \R{1g}.

%% CHANGE-2A bound_cond, init_val

We propose a six-point, two-level Crank-Nicholson type scheme for
solving \R{general_PDE} and \R{1g}-\R{init_val},
\bbbb&&\frac{c_{5}}{2h\tau}\left[u_{m+1,n} - u_{m-1,n} - u_{m+1,n-1} + u_{m-1,n-1} \right] \nnn \\ &&~~~+\frac{c_{4}}{2h^{2}}\left[u_{m+1,n} - 2u_{m,n} +
u_{m-1,n} + u_{m+1,n-1} - 2u_{m,n-1} + u_{m-1,n-1}\right] \nnn\\
&&~~~+\frac{c_{3}}{4h}\left[u_{m+1,n} - u_{m-1,n} +
u_{m+1,n-1} - u_{m-1,n-1}\right] \nnn \\
&&~~~+\frac{c_{2}}{\tau}[u_{m,n} - u_{m,n-1}] +
\frac{c_{1}}{2}\left[u_{m,n} + u_{m,n-1}\right] + c_{0} = 0.\eeee
The above can be conveniently reformulated to a partial difference
equation \bbb&&\left(\frac{c_{5}}{2h\tau} + \frac{c_{4}}{2h^{2}} +
\frac{c_{3}}{4h}\right)u_{m+1,n} +
\left(-\frac{c_{4}}{h^2}+\frac{c_{2}}{\tau}+\frac{c_{1}}{2}\right)u_{m,n}
+ \left(-\frac{c_{5}}{2h\tau} + \frac{c_{4}}{2h^{2}} - \frac{c_{3}}{4h}\right)u_{m-1,n}
\nnn\\ &&~~~= \left(\frac{c_{5}}{2h\tau}-\frac{c_{4}}{2h^{2}} -
\frac{c_{3}}{4h}\right)u_{m+1,n-1} +
\left(\frac{c_{4}}{h^2}+\frac{c_{2}}{\tau}-\frac{c_{1}}{2}\right)u_{m,n-1}
\nnn\\&&~~~~~~~ \left(-\frac{c_{5}}{2h\tau}-\frac{c_{4}}{2h^{2}} +
\frac{c_{3}}{4h}\right)u_{m-1,n-1} +
c_{0}.\label{general_scheme}\eee

Let $w$ be a sufficiently smooth function defined on $\Omega.$ We
define $$ Pw(r,z) = c_{5}\frac{\partial^2
w}{\partial z \partial r} + c_{4}\frac{\partial^2 w}{\partial r^2}(r,z) +
c_{3}\frac{\partial w}{\partial r}(r,z) + c_{2}\frac{\partial
w}{\partial z}(r,z) + c_{1}w(r,z) + c_{0}$$ and \bbb
&&P_{h,\tau}w_{m,n-1/2} = \l(\frac{c_{5}}{2h\tau} + \frac{c_{4}}{2h^{2}} +
\frac{c_{3}}{4h}\r)w_{m+1,n} +
\l(-\frac{c_{4}}{h^2}+\frac{c_{2}}{\tau}+\frac{c_{1}}{2}\r)w_{m,n}
\nnn\\
&&~~+ \l(-\frac{c_{5}}{2h\tau} + \frac{c_{4}}{2h^{2}} - \frac{c_{3}}{4h}\r)w_{m-1,n} +
\l(-\frac{c_{5}}{2h\tau} + \frac{c_{4}}{2h^{2}} + \frac{c_{3}}{4h}\r)w_{m+1,n-1} \nnn\\ &&~~+
\l(-\frac{c_{4}}{h^2}-\frac{c_{2}}{\tau}+\frac{c_{1}}{2}\r)w_{m,n-1} + \l(\frac{c_{5}}{2h\tau}
+ \frac{c_{4}}{2h^{2}} -
\frac{c_{3}}{4h}\r)w_{m-1,n-1} + c_{0}. \label{diff_op}\eee

By expressing the function $w$ at the grid points as Taylor
expansions evaluated at reference point
$\l(r_{m},z_{n-\frac{1}{2}}\r)$ and substituting these into
\R{diff_op}, it can be shown that
\bb{local}\l\|\l(P-P_{h,\tau}\r)w_{m,n-1/2}\r\| =
O\l(h^{2}+h^{2}\tau^{2}+\tau^{2}\r).\ee Let $\sigma=\tau/h$ be
bounded and $\sigma\rightarrow 0$ as $h,~\tau\rightarrow 0.$  Then
\R{local} ensures
\begin{enumerate}\item the consistency of the finite difference scheme
\R{general_scheme};\item a second-order local truncation error of
the scheme \R{general_scheme};\item the numerical stability
depends on the particular coefficients of the differential
equation considered.\end{enumerate}

\begin{Definition}  Consider a homogeneous finite difference scheme
written as a system of linear equations as below:
$$B\mathbf{u}_{n} = C\mathbf{u}_{n-1}$$
or $$\mathbf{u}_{n} = B^{-1}C\mathbf{u}_{n-1}$$ \noindent where
vector $\mathbf{u}_{n} = \left\{u_{k,n}\right\}^{M}_{k=0}$ and the
difference operators $B,C \in \CC^{M\times M}$ are coefficient
matrices.  Let $E = B^{-1}C$. If there exists a constant $K > 0$
independent of $n$, $h$ and $\tau$ such that $\|E^{n}\| \leq K$ for
some norm $\|\cdot\|$, we say that the scheme is stable in the
Lax-Richtmyer sense {\rm\cite{Lax, Morton2}.}
\end{Definition}

%% CHANGE-2A the stability definition is changed

% For this paper, we will employ a related criterion.

Due to the inclusion of a cross-derivative term in our transformed equation, our difference scheme will not be stable in the Lax-Richtmyer sense.  We define a notion of practical stability which holds for a range of propagation step sizes that afford us sufficient resolution in our simulations.

\begin{Definition} \label{OurStability} Let $\rho(B^{-1}C)$ be the spectral radius of kernel matrix
$B^{-1}C$.  If
$$\rho(B^{-1}C) \leq 1$$
at all propagation steps $0<n\tau<T$ with a transverse direction step size
$\epsilon_{0}<h<\epsilon_{1}$ for some $\epsilon_{1} > \epsilon_{0} > 0$, we say that the scheme
$$B\mathbf{u}_{n} = C\mathbf{u}_{n-1}$$ is stable within a parameter range.
\end{Definition}

While this stability condition does not specify a norm for
convergence, it does guarantee that perturbations will not increase
exponentially with $n.$

\begin{Definition}
A matrix $A \in \CC^{n \times n}$ is said to be positive semistable
if every eigenvalue of $A$ has nonnegative real part.
\end{Definition}

\begin{Theorem} \label{GA_semistable_implies_scheme_stable}
Let $A,B,C,G \in \CC^{M\times M}$ be such that $$B = G+A, ~~~C = G-A
.$$
Then the difference scheme defined by
$$B\mathbf{u}_{n} = C\mathbf{u}_{n-1}$$
is stable if and only if $G^{-1}A$ is positive semistable.
\end{Theorem}

\begin{Corollary} \label{dI_A_semistable_implies_scheme_stable}
Let $A,B,C \in \CC^{M\times M},$ $d$ be a positive real number such
that
$$B = dI+A, ~~~C = dI-A.$$
Then the difference scheme defined by $$B\mathbf{u}_{n} =
C\mathbf{u}_{n-1}$$ is stable if and only if $A$ is positive
semistable.
\end{Corollary}

%% CHANGE-2A bound_cond, init_val

Recall \R{1f}. We have the corresponding paraxial Helmholtz
coefficients for the general equations \R{general_PDE} and
\R{general_scheme}, \bb{general_c}c_{5} = 0,~ c_{4}= 1,~ c_{3}= \frac{1}{r},~
c_{2}= -2i\kappa,~ c_{1}= 0,~c_{0} = 0.\ee It follows therefore
that \R{general_scheme}, \R{bound_cond} and \R{init_val} can be
simplified to the following homogeneous paraxial Helmholtz
difference scheme, \bbb&&\hspace{-1cm}-\alpha\left(1 +
\frac{1}{2m}\right)u_{m+1,n} + \left(2 + 2\alpha\right)u_{m,n}
- \alpha\left(1 - \frac{1}{2m}\right)u_{m-1,n} \nnn\\
&&\hspace{-1cm}~~~= \alpha\left(1 + \frac{1}{2m}\right)u_{m+1,n-1}
+ \left(2 - 2\alpha\right)u_{m,n-1} + \alpha\left(1 -
\frac{1}{2m}\right)u_{m-1,n-1},\label{homogeneous_scheme2}\\
&&\hspace{-1cm}u_{m,0}=e^{-hm/\beta_0},\label{initial2}\\
&&\hspace{-1cm}-2\alpha u_{1,n} + \left(2 + 2\alpha\right)u_{0,n}
 = 2\alpha u_{1,n-1} + \left(2 -
 2\alpha\right)u_{0,n-1},\label{boundary2}\\
&&\hspace{-1cm}\left(2 + 2\alpha\right)u_{M,n} - 2\alpha u_{M-1,n}
 = \left(2 - 2\alpha\right)u_{M,n-1} + 2\alpha
 u_{M-1,n-1},\label{boundary3}\eee
where $$\alpha = -\frac{\tau i}{2\kappa h^{2}}.$$

Following the analysis method outlined above, we can express our
scheme in matrix form $$B\mathbf{u}_{n} = C\mathbf{u}_{n-1}$$ where
$B = G+A,$ $C=G-A,$ $G=2I,$ and $A$ is tridiagonal. Investigating
properties of the eigenvalues of $A = \{a_{m,n}\},$ where \bbbb
a_{m,m} &=&
2\alpha,~~~m=0,1,\ldots,M,\\
a_{m,m-1}&=& -\alpha\left(1 -
\frac{1}{2m}\right),~~~m=1,2,\ldots,M-1,\\
a_{M,M-1} &=& -2\alpha,\\
a_{m,m+1}&=& -\alpha\left(1 +
\frac{1}{2m}\right),~~~m=1,2,\ldots,M-1,\\
a_{0,1} &=& -2\alpha,\eeee w are able to show that the eigenvalues
of matrix $A$ are purely imaginary, and thus have nonnegative real
parts. Thus, $A$ is positive semidefinite. By Corollary
\ref{dI_A_semistable_implies_scheme_stable}, we can show the
following.

\vspace{3mm}

%% CHANGE-2A added sentence to theorem.

\begin{Theorem}Let $\kappa$ be a constant. Then the
homogeneous paraxial Helmholtz difference scheme
\R{homogeneous_scheme2}-\R{boundary3} is stable.  Further, there is lower boundary restriction on the step size parameter, $h,$ in this case.
\end{Theorem}

\vspace{8mm}

\noindent{\large\bf 3.~~$z$-Stretching Domain Transformation}
\setcounter{section}{3}\setcounter{equation}{0}

\vspace{1mm}

\noindent One possible way of avoiding the computational
difficulties presented by the discontinuity of $\kappa$ at the
interface is by decomposing the domain into three sections, pre-lens,
lens, and post-lens, and {\em stretching\/} each segment by
one-to-one transformations onto rectangular areas. We would then
be able to use conventional finite difference techniques, such as
that introduced in Section 2, to solve \R{general_PDE} together
with initial-boundary conditions on each segment. The grid stretch
can be achieved either in the direction of electro-magnetic wave
propagation $z,$ or the direction of $r.$ Each of the approaches
have distinct advantages. We will only focus on the former
strategy in this paper. In this case, the numerical solution
computed at the rightmost edge of the pre-lens segment becomes the initial condition of the next segment.

Let $r=r(\xi,\zeta),~z=z(\xi,\zeta)$ be the one-to-one stretching
transformation to be used. Thus, \bbbb&&\frac{\partial u}{\partial
r}=\frac{\partial u}{\partial\xi}\frac{\partial \xi}{\partial
r}+\frac{\partial u}{\partial\zeta}\frac{\partial\zeta}{\partial
r},~~~\frac{\partial u}{\partial z}=\frac{\partial u}{\partial
\xi}\frac{\partial\xi}{\partial z}+\frac{\partial u}{\partial
\zeta}\frac{\partial\zeta}{\partial z},\\&&\frac{\partial^{2}
u}{\partial r^{2}}=\frac{\partial
u}{\partial\xi}\frac{\partial^{2}\xi}{\partial
r^{2}}+\frac{\partial^{2} u}{\partial
{\xi}^{2}}\left(\frac{\partial\xi}{\partial r}\right)^{2} +
2\frac{\partial^2 u}{\partial\xi\partial\zeta}\frac{\partial
\xi}{\partial r}\frac{\partial\zeta}{\partial r}+\frac{\partial
u}{\partial\zeta}\frac{\partial^{2}\zeta}{\partial r^{2}} +
\frac{\partial^{2} u}{\partial {\zeta}^{2}}\left(\frac{\partial
\zeta}{\partial r}\right)^{2}.\eeee A substitution of the above
into \R{1f} yields \bbbb && 2i\kappa\left(\frac{\partial
u}{\partial \xi}\frac{\partial \xi}{\partial z}+\frac{\partial
u}{\partial\zeta}\frac{\partial \zeta}{\partial z}\right) =
\frac{\partial u}{\partial \xi}\frac{\partial^{2}\xi}{\partial
r^{2}}+\frac{\partial^{2} u}{\partial
{\xi}^{2}}\left(\frac{\partial\xi}{\partial r}\right)^{2} +
2\frac{\partial^2 u}{\partial\xi\partial
\zeta}\frac{\partial\xi}{\partial r}\frac{\partial\zeta}{\partial
r}\\&&~~~~ +\frac{\partial
u}{\partial\zeta}\frac{\partial^{2}\zeta}{\partial r^{2}} +
\frac{\partial^{2}
u}{\partial{\zeta}^{2}}\left(\frac{\partial\zeta}{\partial
r}\right)^{2} + \frac{1}{r}\left(\frac{\partial u}{\partial
\xi}\frac{\partial\xi}{\partial r}+\frac{\partial u}{\partial
\zeta}\frac{\partial\zeta}{\partial r}\right),\eeee which can be
regrouped into \bbb&&\left(2i\kappa\frac{\partial \zeta}{\partial
z} - \frac{\partial^{2}\zeta}{\partial r^{2}} -
\frac{1}{r}\frac{\partial\zeta}{\partial r}\right)\frac{\partial
u}{\partial \zeta} = \left(-2i\kappa\frac{\partial \xi}{\partial
z} + \frac{\partial^{2}\xi}{\partial
r^{2}}+\frac{1}{r}\frac{\partial \xi}{\partial
r}\right)\frac{\partial u}{\partial\xi} \nnn\\ &&~~~~+
\left(\frac{\partial \xi}{\partial r}\right)^{2}\frac{\partial^2
u}{\partial {\xi}^2} + 2\left(\frac{\partial \xi}{\partial
r}\frac{\partial \zeta}{\partial r}\right)\frac{\partial^2
u}{\partial\xi\partial\zeta} + \left(\frac{\partial
\zeta}{\partial r}\right)^{2}\frac{\partial^2
u}{\partial{\zeta}^2}.\label{3b}\eee

To map the lens area $\Omega=\l\{0 \leq z \leq Z,~ (z-R)^2 + r^2
\leq R^2,~r \geq 0\r\}$ into a rectangular area
$\tilde{\Omega}=\l\{0 \leq \zeta \leq Z,~0 \leq \xi \leq
R_{1}\r\},$ a natural choice is the following transformation,
$$\xi(r,z) =
r,~~~\zeta(r,z) = \frac{z-R+\sqrt{R^2 - r^2}}{Z-R+\sqrt{R^2 -
r^2}}Z.$$ In this particular case we have \bbbb &&\frac{\partial
\xi}{\partial r} = 1,~~\frac{\partial\xi}{\partial z} = 0,~~
\frac{\partial^2 \xi}{\partial r^2} = 0,\\
&&\frac{\partial \zeta}{\partial r} =
\frac{rZ(z-Z)}{{\rho}^2\sqrt{R^{2} - r^{2}}},~~ \frac{\partial
\zeta}{\partial z} = \frac{Z}{\rho},~\frac{\partial^{2}
\zeta}{\partial r^{2}} = \frac{Z (z-Z)\l(\rho R^2+2 r^{2}
\sqrt{R^{2}-r^{2}}\r)}{{\rho}^3\l(R^{2}-r^{2}\r)^{3/2}},\eeee where
$\rho=Z-R+\sqrt{R^2-r^2}.$ We define functions \bbbb\phi(\xi,\zeta)
&:=& \frac{\partial\zeta}{\partial r}\l(r\l(\xi,\zeta\r),
z\l(\xi,\zeta\r)\r) = \frac{\l(\zeta-Z\r)\xi}{\sqrt{R^{2}-\xi^{2}}
\l[Z-\l(R-\sqrt{R^{2}-\xi^{2}}\r)\r]}, \\
\psi(\xi,\zeta) &:=& \frac{\partial^{2}\zeta}{\partial
r^{2}}\l(r\l(\xi,\zeta\r),z\l(\xi,\zeta\r)\r) =
\frac{\l(\zeta-Z\r)\l[R^{3}-\sqrt{R^{2}-\xi^{2}}
\l(R^{2}+2\xi^{2}\r)-Z\r]}{\l(R^{2}-\xi^{2}\r)^{\frac{3}{2}}
\l[Z-\l(R-\sqrt{R^{2}-\xi^{2}}\r)\r]^{2}},\\
\theta(\xi,\zeta) &:=& \frac{\partial\zeta}{\partial
z}\l(r\l(\xi,\zeta\r),z\l(\xi,\zeta\r)\r) =
\frac{Z}{Z-\left(R-\sqrt{R^{2}-\xi^{2}}\right)}. \eeee Subsequently,
\R{3b} can be simplified to \bb{ZStretchedHelmholtz} \l(2i\kappa
\theta - \psi - \frac{1}{\xi}\phi\r)\frac{\partial u}{\partial\zeta}
= \frac{1}{\xi} \frac{\partial u}{\partial \xi} + \frac{\partial^2
u}{\partial {\xi}^2} +
2\phi\frac{\partial^2u}{\partial\xi\partial\zeta} +
\phi^2\frac{\partial^2 u}{\partial {\zeta}^2}.\ee  It can be
demonstrated that at every point in the transformed lens segment, we
have
$$\phi^{2}\l|\frac{\partial^2 u}{\partial {\zeta}^2}\r| \leq \frac{R_{1}^{2}Z}{\left(R^{2} - R_{1}^{2}\right)
\left[Z-\left(R-\sqrt{R^{2}-R_{1}^{2}}\right)\right]}\l|\frac{\partial^2 u}{\partial z^2}\r|
\approx 0.$$

Thus for a typical lens where the slowly varying
envelope approximation is applicable, \R{ZStretchedHelmholtz} can
be simplified to a more appropriate form for computations within
the transformed lens area, \bb{SimplifiedZStretchedHelmholtz}
\l(2i\kappa
\theta - \psi - \frac{1}{\xi}\phi\r)\frac{\partial u}{\partial\zeta} = \frac{1}{\xi}
\frac{\partial u}{\partial \xi} + \frac{\partial^2 u}{\partial
{\xi}^2} + 2\phi\frac{\partial^2u}{\partial\xi\partial\zeta}.\ee

\begin{center}
\epsfig{file=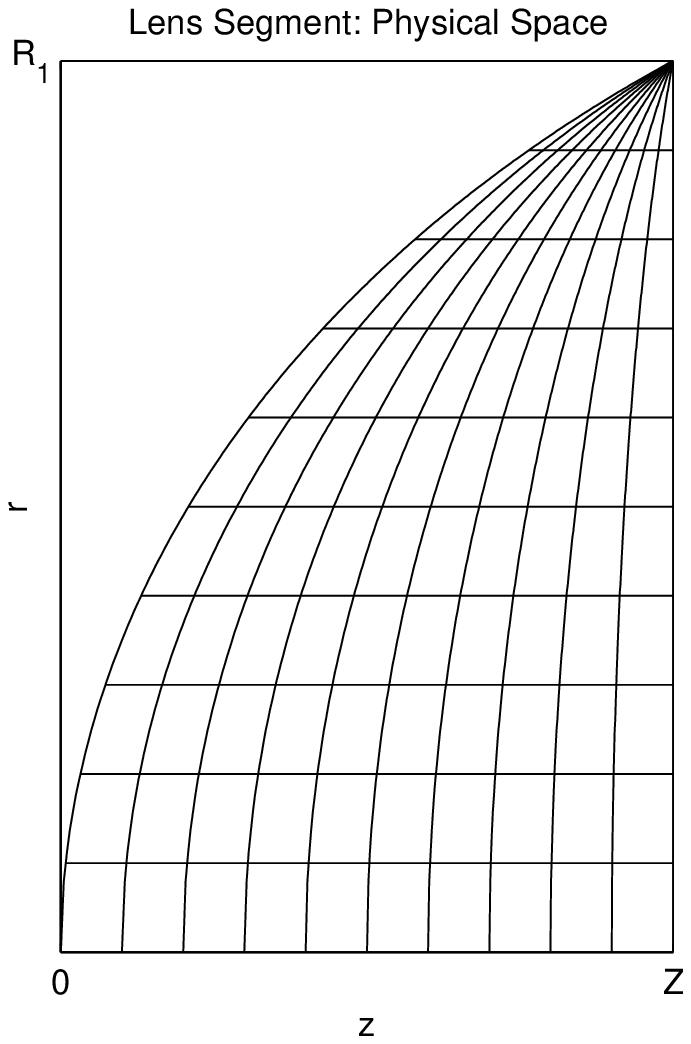,width=2.5in,height=2.35in}
\epsfig{file=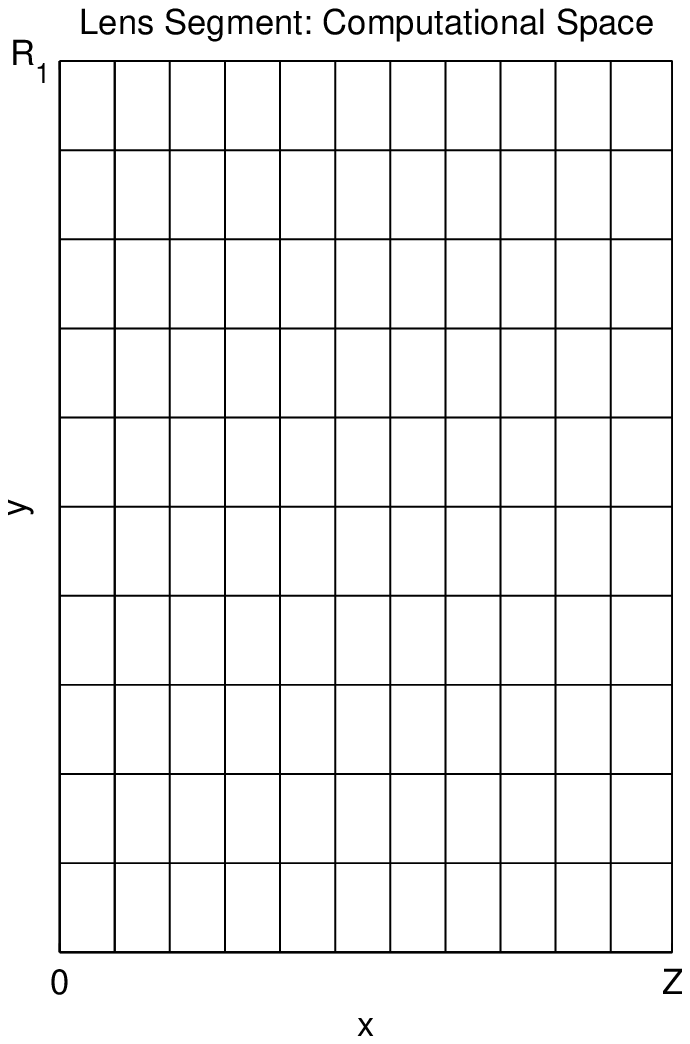,width=2.5in,height=2.35in}

\parbox[t]{13cm}{\small{\sc Figure 3.1}. ~LEFT: An in-lens
domain before a $z$-stretching. RIGHT: The in-lens domain after a
$z$-stretching.}
\end{center}

In Figure 3.1, we illustrate the $z$-stretching by showing an
example of the grid within the lens before and after the
stretching. It is evident that the transformation enables the
six-point scheme to work accurately and efficiently. Instead of
\R{general_c}, now we have the following set of adjusted
coefficients for the in-lens segment difference scheme,
$$c_{5} = 2\phi,~ c_{4} = 1,~ c_{3} = \frac{1}{\xi},~ c_{2} =
-2i\kappa\theta +
\psi + \frac{1}{\xi}\phi,~c_{1} = 0,~c_{0} = 0.$$ Subsequently,
\R{homogeneous_scheme2} can be modified to
\bbb&&\hspace{-1.2cm}-\gamma\l[\frac{2\phi}{h}+\alpha\l(1 +
\frac{1}{2m}\r)\r]u_{m+1,n} + (2+2\alpha\gamma)u_{m,n} -
\gamma\l[-\frac{2\phi}{h}+\alpha\left(1 -
\frac{1}{2m}\right)\r]u_{m-1,n}\nnn\\
&&\hspace{-1.2cm}~~~ = \gamma\l[-\frac{2\phi}{h}+\alpha\left(1 +
\frac{1}{2m}\right)\r]u_{m+1,n-1} + (2-2\alpha\gamma)u_{m,n-1} \nnn\\ ~~~~&&+
\gamma\l[\frac{2\phi}{h}+\alpha\left(1 -
\frac{1}{2m}\right)\r]u_{m-1,n-1},\label{scheme_in_lens}\eee where
$$\alpha = \frac{\tau}{h^{2}} ~\mbox{ and }~
\gamma=\gamma(\xi,\zeta) = \left(2i\kappa\theta-\psi-\frac{1}{\xi}\phi\r)^{-1}.$$

In this method, the Gaussian beam input equation is evaluated at the lens surface,
and becomes the initial solution of the simulation at the edge of the lens segment.  No computation is necessary in the pre-lens segment.  The solution at the right edge
of the lens segment becomes the initial solution of the post-lens segment, which we
simulate using the homogeneous scheme described earlier.

To determine the boundary conditions applicable within the lens segment, note that
$$\frac{\partial u}{\partial r}(\xi,\zeta) = \frac{\partial
u}{\partial \xi}(\xi,\zeta)\frac{\partial \xi}{\partial
r}(\xi,\zeta) + \frac{\partial u}{\partial
\zeta}(\xi,\zeta)\frac{\partial \zeta}{\partial r}(\xi,\zeta),$$
thus \bb{boundary4}\frac{\partial u}{\partial\xi}(0,\zeta) =
\frac{\partial u}{\partial\xi}(R_{1},\zeta) = 0,~~~0<\zeta<Z.\ee It is
interesting to note that for a lens that tapers to a point at the
top, the geometric interpretation of this new boundary condition is
that the single point $(R_{1},Z)$ has been stretched into the upper
edge of our transformed rectangular domain, i.e. the upper boundary
in the computational space corresponds to the single point
$(R_{1},Z)$ in the physical space.

To demonstrate stability, we use the matrix analysis method introduced in the previous section.  For our scheme $B\mathbf{u}_{n} = C\mathbf{u}_{n-1}$, we have $B = G+A,$ $C=G-A$ with matrix $G = \{g_{m,n}\}$ and $A = \{a_{m,n}\},$ where
\bbbb g_{m,m} &=& 2,~~~m=0,1,\ldots,M,\\
g_{m,m-1}&=& \frac{2\gamma\phi}{h},~~~m=1,2,\ldots,M-1,\\
g_{M,M-1} &=& 0,\\
g_{m,m+1}&=& \frac{-2\gamma\phi}{h},~~~m=1,2,\ldots,M-1,\\
g_{0,1} &=& 0,\eeee
\bbbb a_{m,m} &=&
2\alpha\gamma,~~~m=0,1,\ldots,M,\\
a_{m,m-1}&=& -\alpha\gamma\left(1 -
\frac{1}{2m}\right),~~~m=1,2,\ldots,M-1,\\
a_{M,M-1} &=& -2\alpha\gamma,\\
a_{m,m+1}&=& -\alpha\gamma\left(1 +
\frac{1}{2m}\right),~~~m=1,2,\ldots,M-1,\\
a_{0,1} &=& -2\alpha\gamma.\eeee

%% CHANGE-2A typo:  $G^{-1}A$ is positive semidefinite -> $G^{-1}A$ is positive semistable
%% CHANGE-2A typo:  $A$ is positive semidefinite -> $A$ is positive semistable
%% CHANGE-2A replaced paragraph, added paragraphs following theorem and table

%By examining properties of $\phi\l(\xi,\zeta\r)$ and
%$\gamma\l(\xi,\zeta\r)$ for our particular transformation, we are
%able to conclude that within useful parameter ranges, $G+G^{*}$ is
%positive definite and $A$ is positive semistable, thus $G^{-1}A$
%is positive semistable. Based on Theorem
%(\ref{GA_semistable_implies_scheme_stable}), we can further prove
%the following.

Examining the real part of function $\gamma\l(\xi,\zeta\r)$ resulting from our chosen transformation, we see that matrix $A$ is positive semistable.  This property will still hold if the transformation is adapted to other convex lens shapes.  Further, if we have a transverse step size $h$ such that
$$h > 2\l|\gamma\l(\xi,\zeta\r)\phi\l(\xi,\zeta\r)\r| ~\text{ for }~ 0 \leq \xi \leq R_{1}, ~0 \leq \zeta \leq Z$$
\noindent or equivalently, if the number of grid points in the $\xi$ direction, $M,$ is such that
$$M < \frac{R_{1}}{h_{min}}$$
\noindent where $$h_{min} = 2\max_{\xi,\zeta}\l|\gamma\l(\xi,\zeta\r)\phi\l(\xi,\zeta\r)\r|$$
\noindent then matrix $G+G^{*}$ is positive definite.  Then based on Theorem
(\ref{GA_semistable_implies_scheme_stable}), we can prove the following.

\vspace{3mm}

\no{\bf Theorem 3.1.} {\em Let $\kappa$ be discontinuous as given by
\R{1g} and
$$h > 2\l|\gamma\l(\xi,\zeta\r)\phi\l(\xi,\zeta\r)\r| ~\text{ for }~ 0 \leq \xi \leq R_{1}, ~0 \leq \zeta \leq Z.$$  Then the difference scheme
\R{scheme_in_lens}-\R{boundary4} is stable on $z$-stretched
domains.}

\vspace{3mm}

For the parameter values utilized in our simulations
\begin{equation}\label{our_params} \kappa = 9.97543 \times 10^{3}, ~~ R = 1.969, ~~ Z = 0.7643, ~~ R_{1} = 1.5574
\end{equation}

\noindent we need $M \leq 1.2092 \times 10^{4}.$  We list a sampling of maximum grid points $M$ for other parameter values.

\begin{center}
\begin{tabular}{|c|ccc|}
    \hline
$Z$ &  $k=8000$    &   $k=10000$  & $k=12000$ \\
    \hline
0.1   & 7439 & 9365 & 11279  \\
0.3   & 5642 & 7046 & 8848  \\
0.5   & 4059 & 5064 & 6068 \\
0.7   & 2496 & 3104 & 3711  \\
0.9   & 1033 & 1257 & 1479  \\
    \hline
\end{tabular}

\parbox[t]{13cm}{\small{\sc Table 3.1}. ~Maximum grid points in the transverse direction with $R=1.$}
\end{center}

\vspace{4mm}

\noindent{\large\bf 4.~~Simulation Results and Observations}
\setcounter{section}{4}\setcounter{equation}{0}

\vspace{1mm}

%% CHANGE-2A typo sack -> sake, added C++

\no All simulated results were implemented on dual processor {\sc
Dell} workstations with at least double precision. {\sc MatLab}, {\sc Fortran} and {\sc C++} programming languages were utilized. Dimensionless
models are used throughout the computations. For the sake of
simplicity, we do not tend to re-scale numerical solutions back to
their original physical dimensions in simulations.

%% CHANGE-2A replace paragraph

%In the following numerical experiments, we choose $r_0=1.56$ and
%$Z=1.53$ for the initial-boundary value problem
%\R{1f}-\R{bound_cond}, while replacing \R{initial} by a normalized
%complex function \bb{initial5}u(r,0)=e^{i\kappa r},~~~0\leq r\leq
%r_0,\ee which is illustrated in Figure 4.1. An auxiliary smoothing
%function \cite{Sheng2} is used with a bandwidth
%$\epsilon_0=10^{-5}.$ We further select $h=\tau=10^{-3}$ in the
%six-point scheme used in the computational solution space after a
%designated $z$-stretching.

In the following numerical experiments we use parameters (\ref{our_params}) listed in the previous section.  We further select $h = R_{1}/M$ where $M = 5 \times 10^{3}$ and $\tau = Z/N$ where $N = 1.6 \times 10^{4}$ in the six-point scheme used in the computational solution space after a designated $z$-stretching.

\vfill\pagebreak

\begin{center}
\epsfig{file=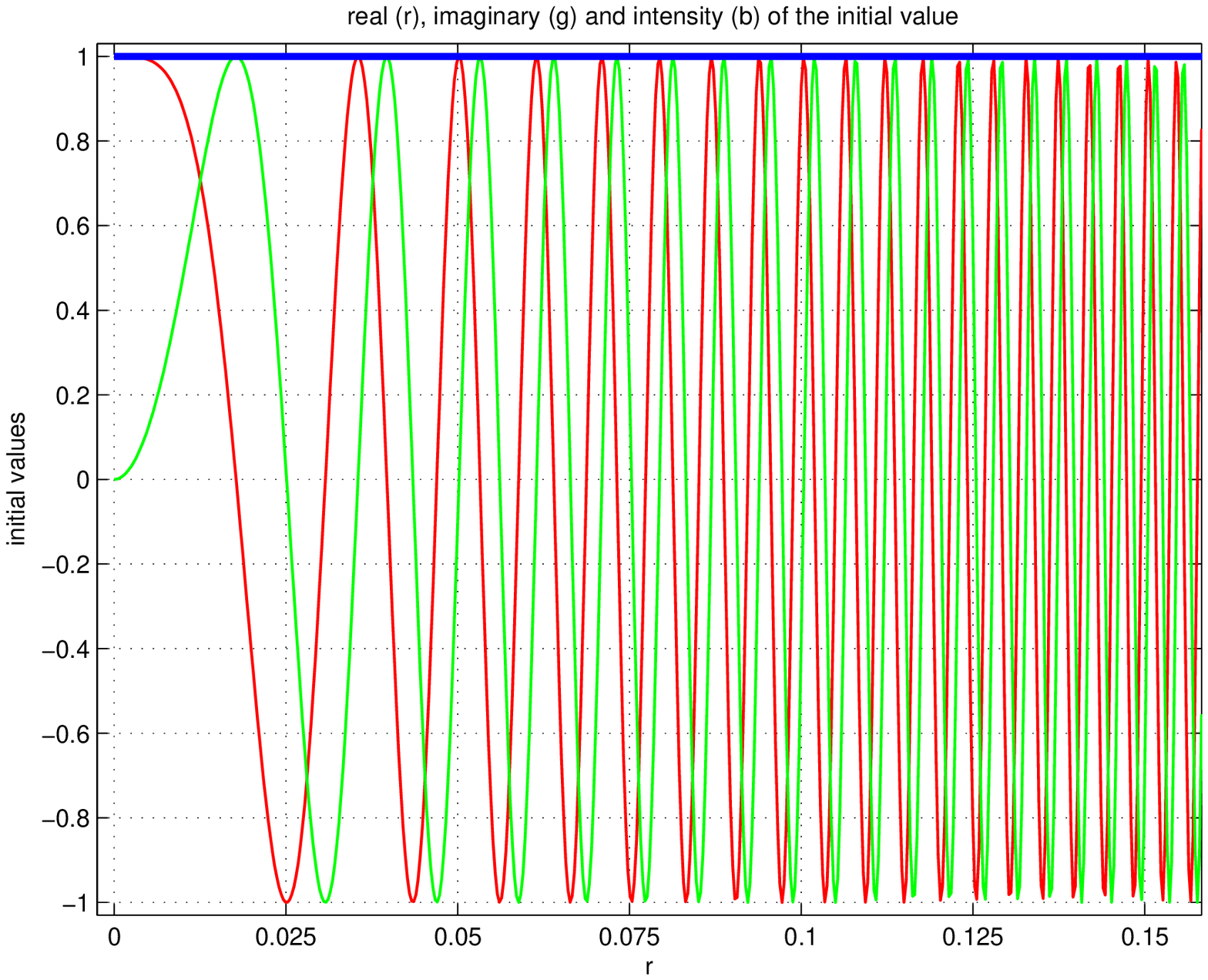,width=4.8in,height=2.4in}

\parbox[t]{13cm}{\small{\sc Figure 4.1}. ~Normalized initial value function $u(r,0).$
Red curve is for the real part and green curve is for the
imaginary part of the function. Highly oscillatory features of the
function is clear.}
\end{center}

We show the real part of the simulated solution,
$a(z)=\mbox{real}\{u(0,z)\},$ in Figures 4.2-3. We may observe that
while the function value of $a$ is relatively stable before the
focusing point $z_f\approx 0.94778,$ it increases dramatically as
$z\rightarrow z_f.$ This can be viewed more precisely in the
enlarged picture of Figure 4.3.

\begin{center}
\epsfig{file=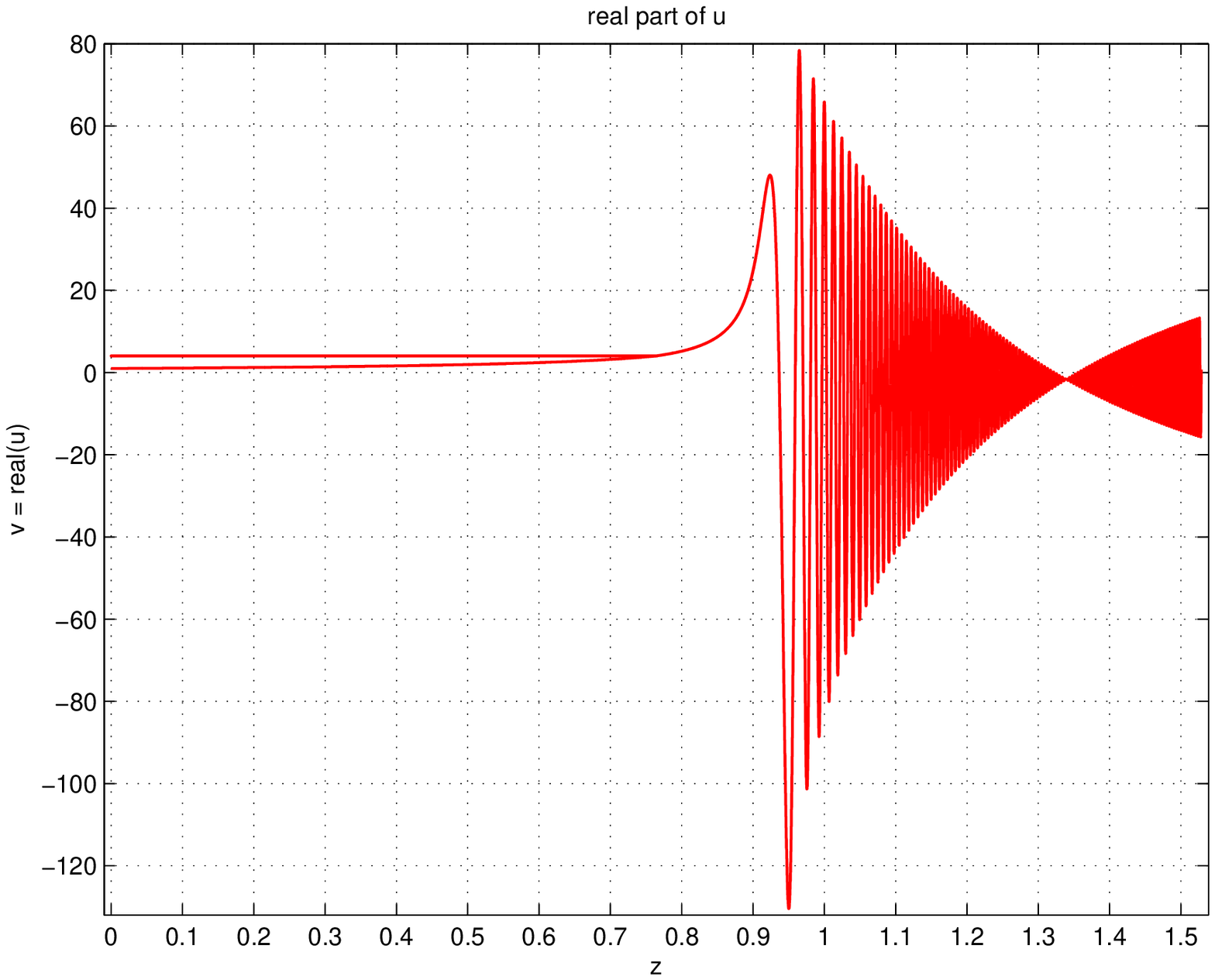,width=4.8in,height=2.4in}

\parbox[t]{13cm}{\small{\sc Figure 4.2}. ~Real part of the simulated
solution at the center point $r=0.$ The numerical solution
increases rapidly as $z$ approaches the focusing location. Then
the simulated oscillatory wave diffuses after the focusing point.
}
\end{center}

\vfill\pagebreak

\begin{center}
\epsfig{file=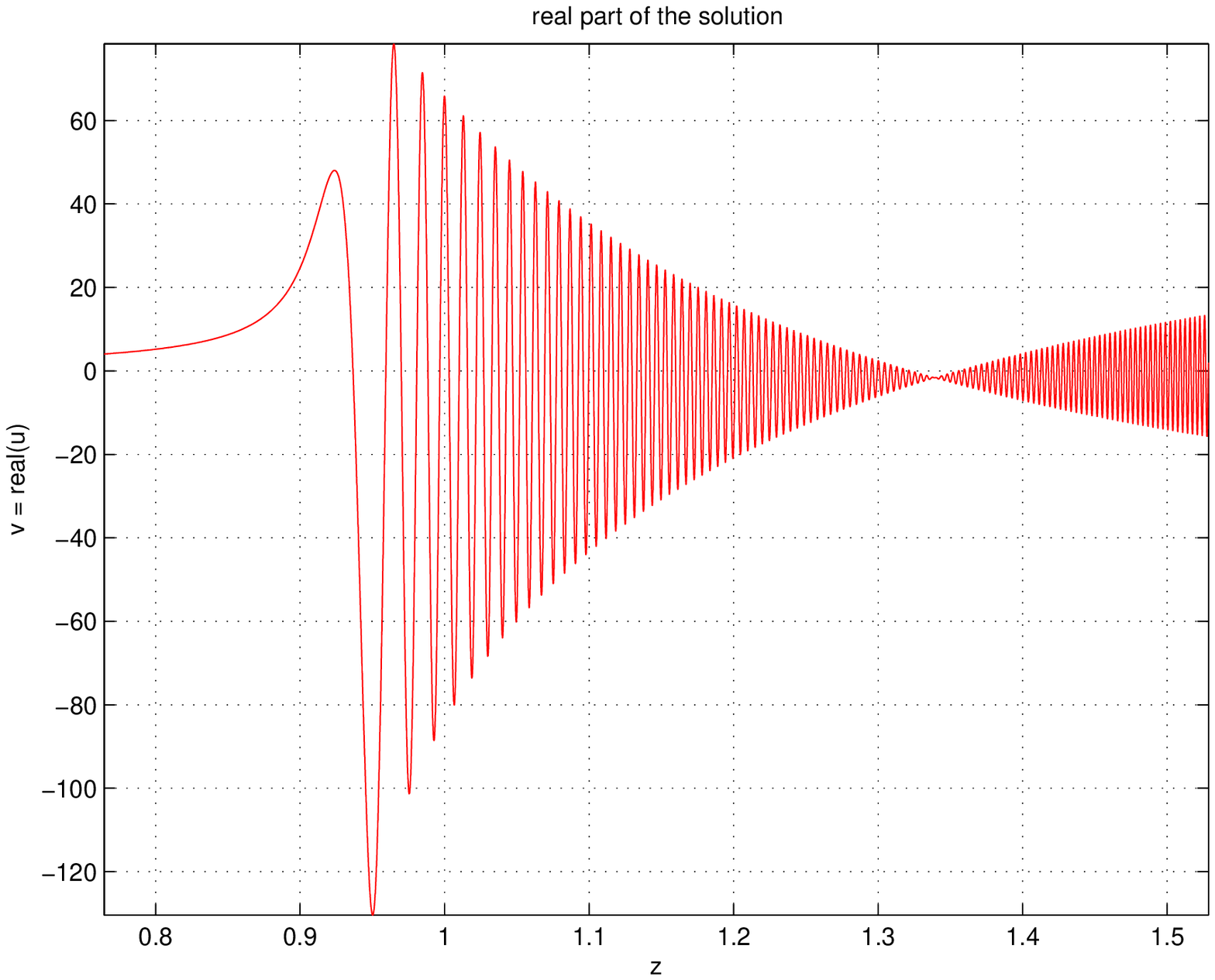,width=4.8in,height=2.4in}

\parbox[t]{13cm}{\small{\sc Figure 4.3}. ~More detailed real part of the simulated
solution near the focus point, $r=0.$ Same conditions as in Figure
4.2 are used. }
\end{center}

%% CHANGE-2A $z_f\approx 2.7431$
Figures 4.4-5 are devoted to the imaginary part of the numerical
solution, $b(z)=\mbox{imaginary}\{u(0,z)\}.$ Similar to the real
part, $b$ is relatively stable before the focusing point $z_f\approx
2.7431$ and is highly oscillatory as $z\rightarrow z_f.$ The
phenomenon can be viewed more clearly in the enlarged picture of
Figure 4.5.

\begin{center}
\epsfig{file=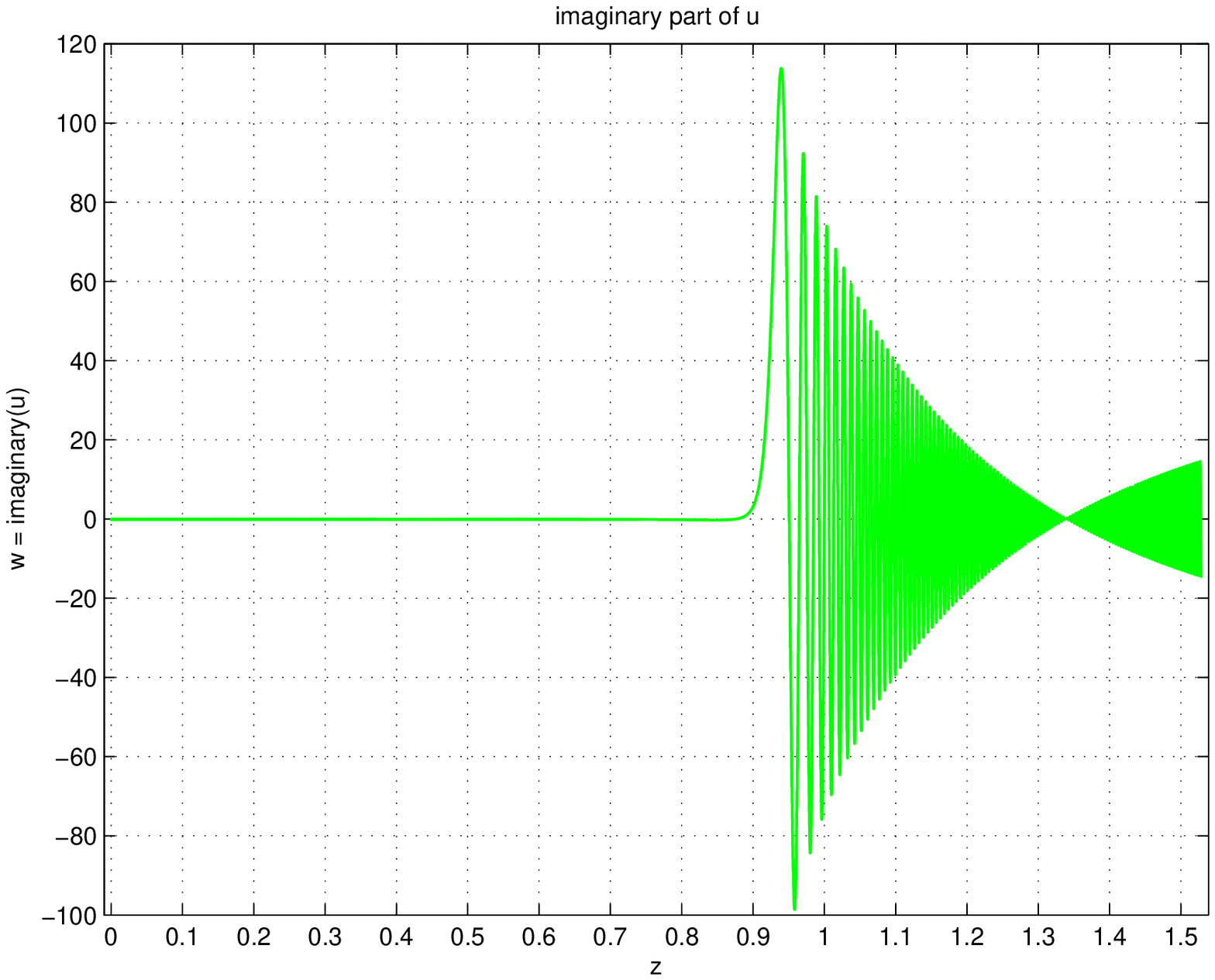,width=4.8in,height=2.4in}

\parbox[t]{13cm}{\small{\sc Figure 4.4}. ~Imaginary part of the simulated
solution at the center point $r=0.$ The numerical solution
increases rapidly as $z$ approaches the focusing location. Then
the simulated oscillatory wave diffuses after the focusing point.}
\end{center}

\vfill\pagebreak

\begin{center}
\epsfig{file=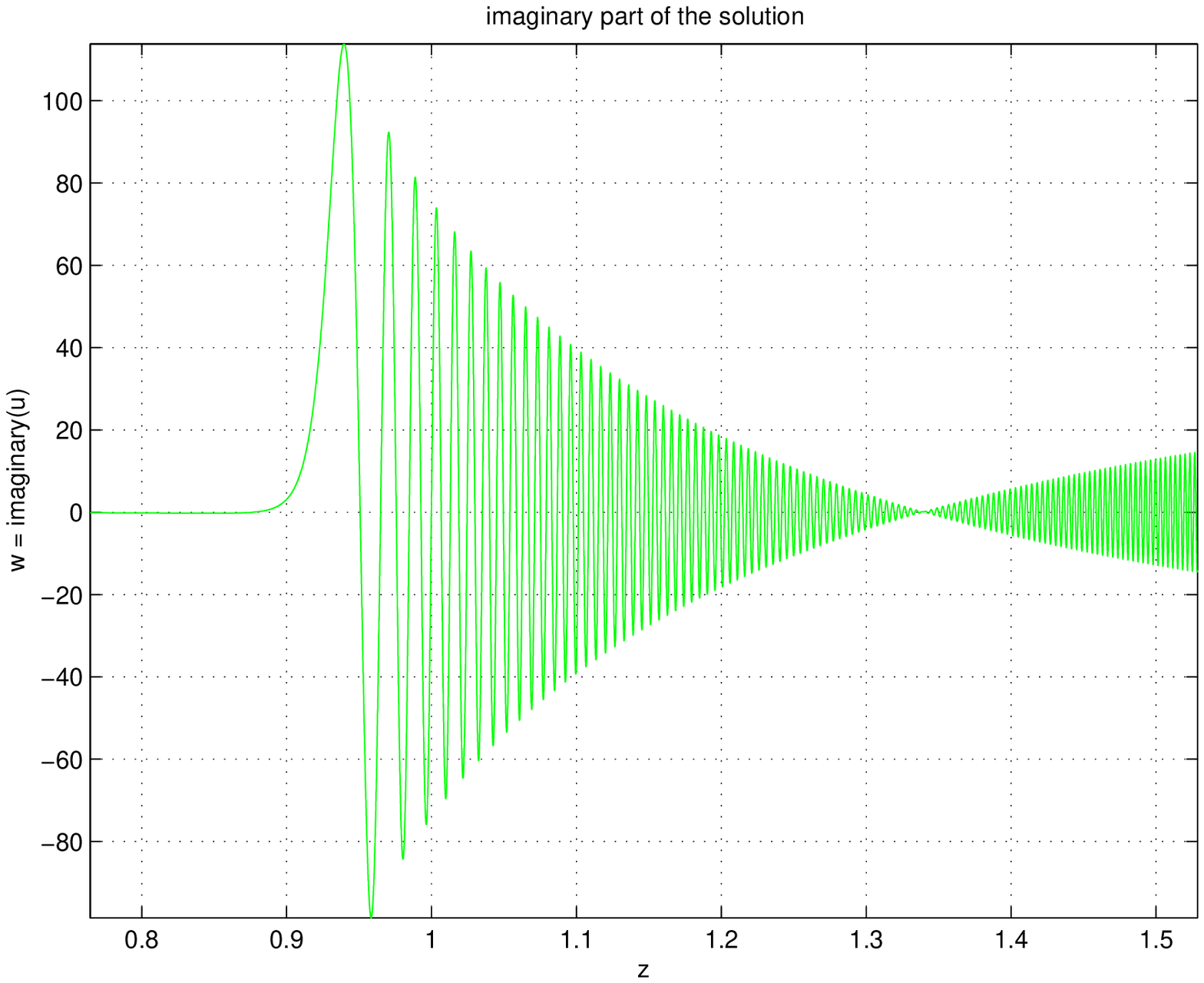,width=4.8in,height=2.4in}

\parbox[t]{13cm}{\small{\sc Figure 4.5}. ~More detailed imaginary part of the simulated
solution near the focus point, $r=0.$ Same conditions as in Figure
4.4 are used. }
\end{center}

Define the {\em numerical intensity function\/} as
\bb{intensity0}{\cal
T}(r,z)=\sqrt{\mbox{real}^2[u(r,z)]+\mbox{imag}^2[u(r,z)]}\geq
0,~~~0\leq r\leq R_{1},~0\leq z\leq Z.\ee In Figures 4.6-7 we plot
this intensity function against the propagation direction $z$ as $r$
being chosen as zero. It is interesting to find that the intensity
increases rapidly as $z\rightarrow z_f.$ The observation is
consistent with our previous results.

\begin{center}
\epsfig{file=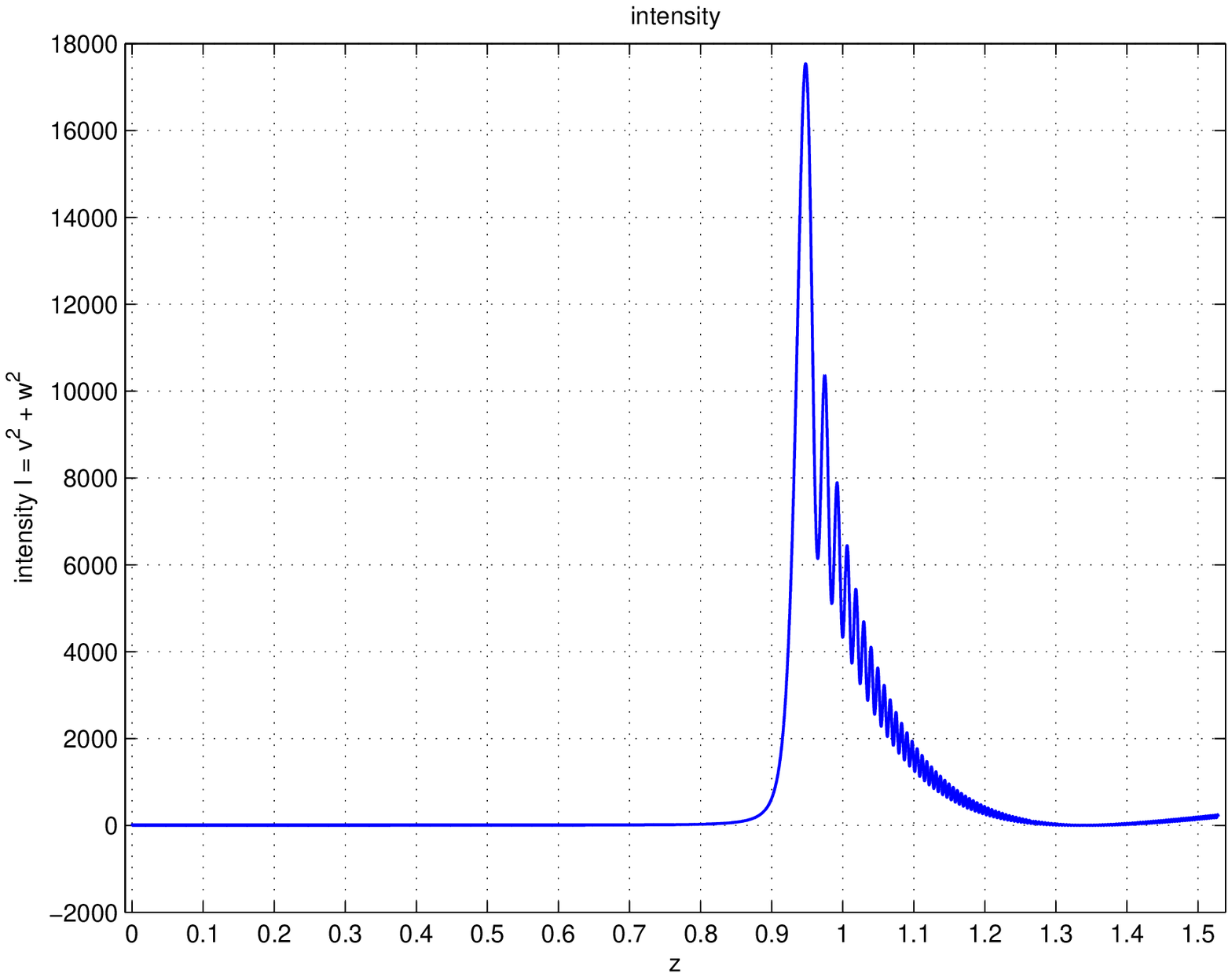,width=4.8in,height=2.4in}

\parbox[t]{13cm}{\small{\sc Figure 4.6}. ~Numerical intensity function of the
simulated solution at the center point $r=0.$ The intensity
increases rapidly as $z$ approaches the focusing location. Then
the intensity value diffuses out after the focusing point. }
\end{center}

\vfill\pagebreak

\begin{center}
\epsfig{file=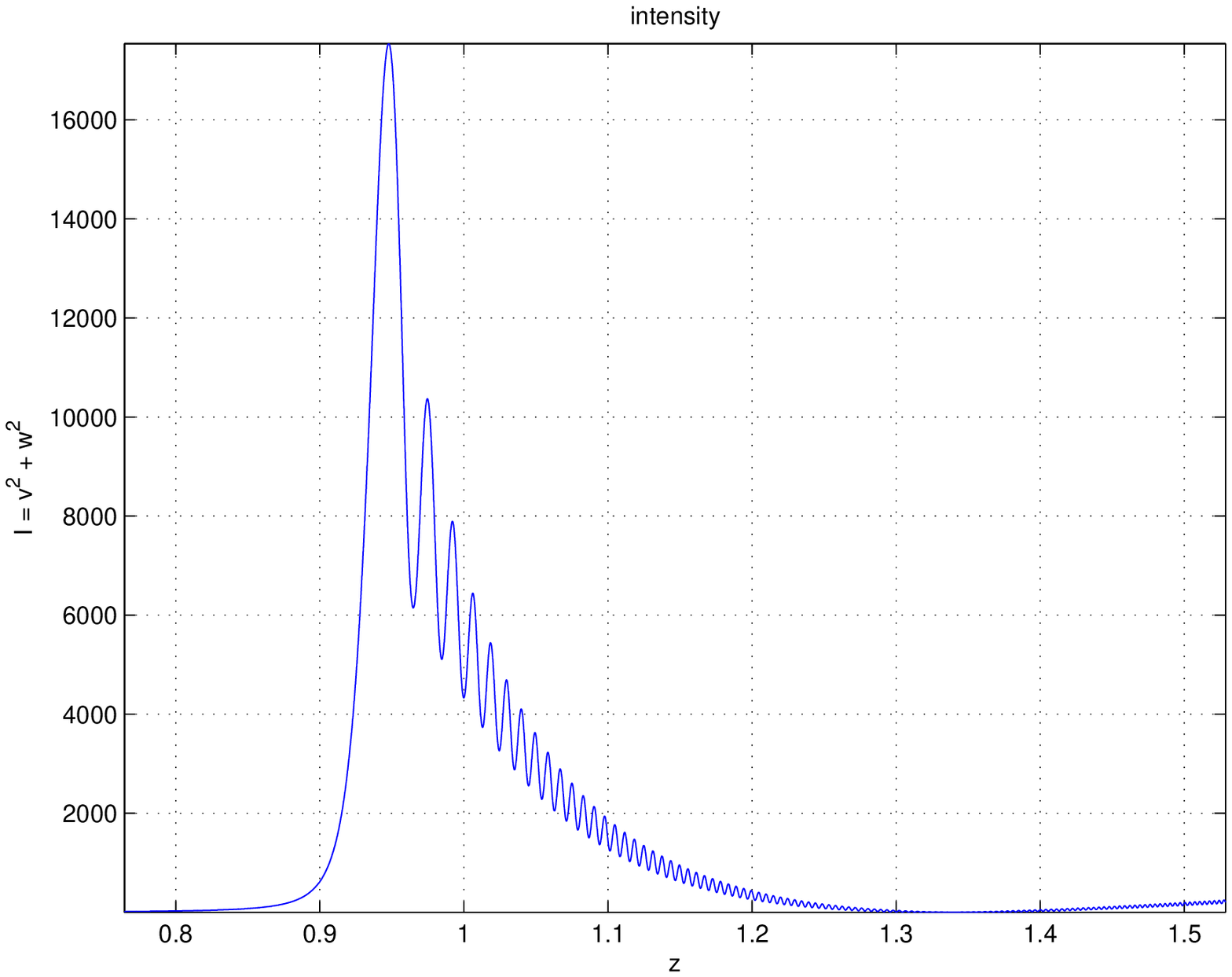,width=4.8in,height=2.4in}

\parbox[t]{13cm}{\small{\sc Figure 4.7}. ~More detailed numerical intensity function
of the simulated solution near the focus point $(r=0).$ Same
conditions as in Figure 4.6 are used.  }
\end{center}

As a comparison, we further plot the numerical intensity function
at the locations near and at the focusing point in Figure 4.8. It
is observed that the numerical estimate of the intensity
oscillates rapidly in the $r$-direction. The intensity increases
sharply near the center point of the lens, $r=0,$ while $z$
approaches the focusing point location. The simulated wave
profiles well match the experimental results. The algorithms can
be used to provide reference values for further explorations.

%% CHANGE-2A remove this figure and caption

%\begin{center}
%\epsfig{file=proj_z1.eps,width=4.8in,height=2.45in}

%\parbox[t]{13cm}{\small{\sc Figure 4.8}. ~Highly oscillatory numerical intensity function of
%the simulated solution at $z=0.7643$ at the right edge of the lens. The intensity builds up near the center of the lens at
%$r=0.$ A domain of $-R_{1}\leq r\leq R_{1}$ is used. }
%\end{center}

%\vfill\pagebreak

\begin{center}
\epsfig{file=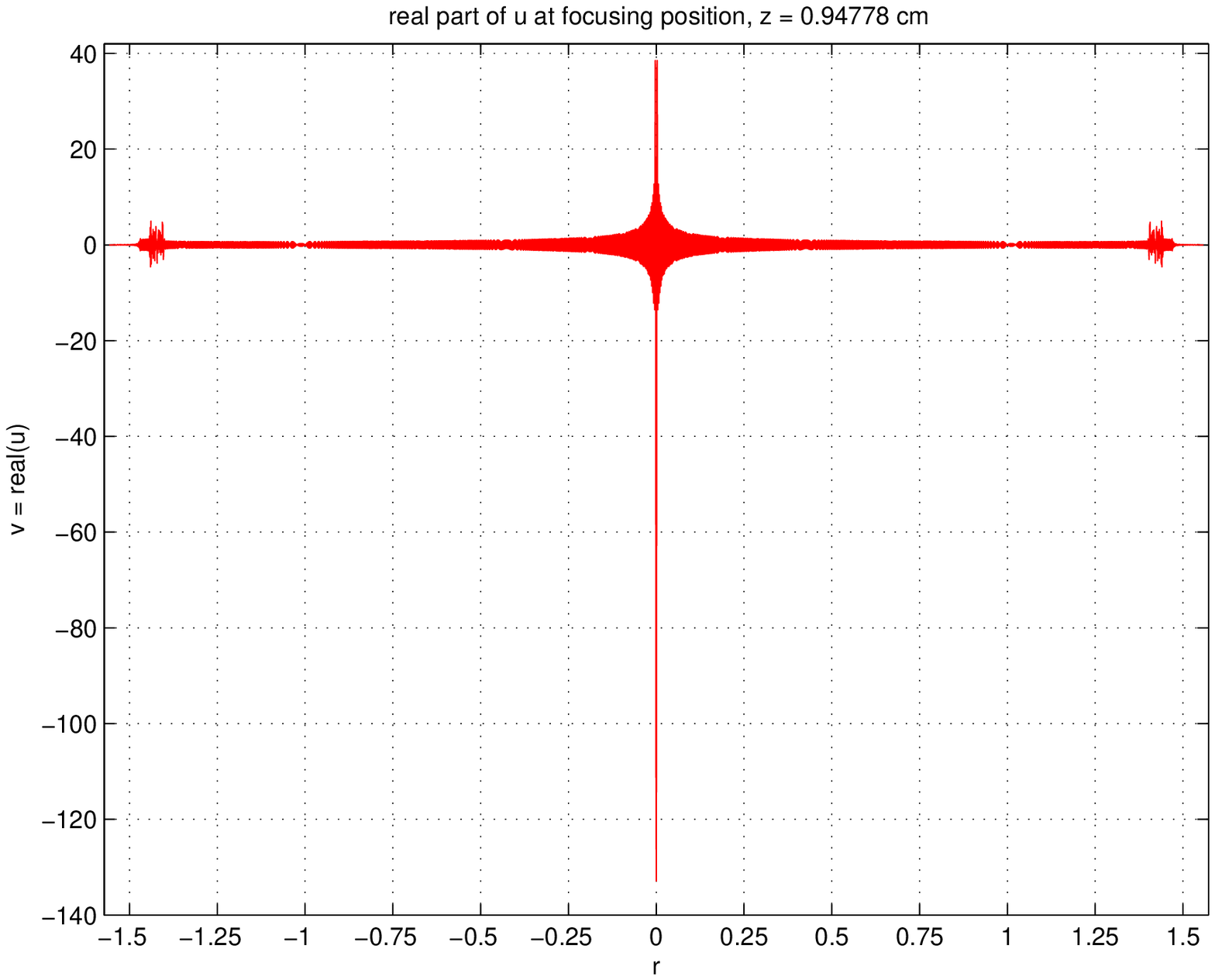,width=4.8in,height=2.45in}

%% CHANGE-2A focusing location $z=0.94778.$ -> focusing location $z=2.7431.$

\parbox[t]{13cm}{\small{\sc Figure 4.8}. ~(ask authors for the image) Highly oscillatory numerical intensity function of
the simulated solution at the focusing location $z=2.7431.$ The
computed intensity increases exponentially near the center, as
compared with much lower profiles away from the focusing area. }
\end{center}

\vspace{6mm}

\noindent{\large\bf 5.~~Conclusions}
\setcounter{section}{5}\setcounter{equation}{0}

\vspace{1mm}

%% CHANGE-2A we don't map pre-lens, removed respectively

In this discussion, we have employed a {\em z-stretching domain
transformation\/} to map the lens and post-lens domains
into convenient rectangular shapes, where we can
utilize well established finite difference methods on uniform grids.
The resulting method provides a useful approximation technique that
can be efficiently implemented with less than 50 lines of Matlab code
in the simulation loops.

A more powerful strategy may be the use of optimally combined $z$
and $r$ stretching transformations to increase computational
resolution and accuracy of the numerical solution in critical local
regions. Magnifying a particular subregion in the transformed
coordinate space is computationally equivalent to increasing the
refinement of the grid in that region using techniques such as
adaptive mesh refinement \cite{Berger1,Sheng}. Alternatively, domain
transformation could facilitate the use of established adaptive grid
techniques, since applications of standard mesh refinement or moving
mesh technics are straightforward on rectangular domains
\cite{Kim,Sheng}. Further, hyperbolic smoothness maps have also been
proved to be extremely useful auxiliary tools to consider in
practical optical beam computations \cite{Sheng2}. Detailed
discussions will be given in our forthcoming reports.

\end{document}